\documentclass[preprint,12pt]{elsarticle}

\usepackage{amssymb, palatino, amsmath, float, caption, subcaption, graphicx, tikz, xcolor}

\journal{journal}

\begin{document}

\begin{frontmatter}



\title{Eco-Evolutionary Dynamics of Bimatrix Games}


\affiliation[inst1]{organization={Department of Mathematics},
            addressline={Dartmouth College}, 
            city={Hanover},
            postcode={03755}, 
            state={NH},
            country={USA}}
            
\affiliation[inst2]{organization={Department of Biomedical Data Science},
            addressline={Geisel School of Medicine at Dartmouth}, 
            city={Lebanon},
            postcode={03756}, 
            state={NH},
            country={USA}}
\author[inst1]{Longmei Shu}

\author[inst1,inst2]{Feng Fu}



\begin{abstract}
Feedbacks between strategies and the environment are common in social-ecological, evolutionary-ecological, and even psychological-economic systems. Utilizing common resources is always a dilemma for community members, like tragedy of the commons. Here we consider replicator dynamics with feedback-evolving games, where the payoffs switch between two different matrices. Although each payoff matrix on its own represents an environment where cooperators and defectors can't coexist stably, we show that it's possible to design appropriate switching control laws and achieve persistent oscillations of strategy abundance. This result should help guide the widespread problem of population state control in microbial experiments and other social problems with eco-evolutionary feedback loops.
\end{abstract}




\end{frontmatter}


\section{Introduction}

Game theory is based on the principle that individuals make rational decisions regarding their choice of actions given suitable incentives~\cite{hauert2005game}. In practice, the incentives are represented as strategy-dependent payoffs. Evolutionary game theory extends game-theoretic principles to model dynamic changes in the frequency of strategies~\cite{nowak2006evolutionary,roca2009evolutionary}. Replicator dynamics is one commonly used framework for such models~\cite{cressman2014replicator}. In replicator dynamics, the frequency of strategists changes over time as a function of the social makeup of the community. The frequency of a strategy increases when the fitness of those who adopt it is greater than the average fitness of the population.

But individual actions do not only modify the social makeup, they very often modify the environment as well~\cite{perc2010coevolutionary}. The strategies that individuals employ impact the environment through time, and the state of the environment in turn influences the payoffs of the game. This environment-dependent feedback occurs across scales from microbes \cite{sanchez2013feedback, wang2020steering, estrela2019environmentally, rand2021geometry} to humans \cite{weitz2016oscillating, hilbe2018evolution, chen2018punishment, shao2019evolutionary, hauert2019asymmetric, tilman2020evolutionary} in public goods game and the tragedy of the commons~\cite{hardin1968tragedy}. The tragedy of the commons is a fundamental problem that is intensively studied~\cite{wu2018coevolutionary,perez2022cooperation,wang2022replicator,wang2022decentralized}. A well-known example is the danger of overfishing. Fishermen are motivated to catch the maximum amount of fish because restraint could only work if all others were behaving similarly. Otherwise, fish are driven to extinction which is the worst scenario for everyone. Similarly, in overgrazing of common pasture lands, an individual's short-term benefit seems to be in conflict with the long-term interest of a larger population. The common feature of these cases is human activity influences the actual state of resources which has a negative feedback for not only those who degrade the environment but also for the whole community. 

The feedback loop between the environment and individual behavioral strategies also appears in microbial populations~\cite{tarnita2017ecology}, including microbes, bacteria, and viruses. Among microbes, feedback may arise due to fixation of inorganic nutrients given depleted organic nutrient availability, the production of extracellular nutrient-scavenging enzymes like siderophores or enzymes like invertase that hydrolyze diffusible products, and the release of extracellular antibiotic compounds. The incentive for public goods production changes as the production influences the environmental state.

Such joint influence occurs in human societies~\cite{cardillo2020critical}, for example, when individuals decide to vaccinate or not~\cite{chen2019imperfect}. Decisions not to vaccinate have been linked most recently to outbreaks of otherwise-preventable childhood infectious diseases in northern California. These outbreaks modify the subsequent incentives for vaccination. Such coupled feedback also arises in public goods dilemmas involving water~\cite{ostrom1993coping} or other resource use such  as antibiotic use~\cite{chen2018social}. In a period of replete resources, there is less incentive for restraint. However, overuse in times of replete resource availability can lead to depletion of the resource and changes in incentives.

A new framework of replicator dynamics with feedback-evolving games has been proposed \cite{weitz2016oscillating} to characterize the phenomenon that the environment and individual behavior co-evolve in many social-ecological and psychological-economic systems~\cite{cao2021eco}. The environmental feedback can result in oscillating dynamics for both the environment quality and strategy states, for example, multiple waves of infections as seen in the COVID-19 pandemic~\cite{glaubitz2020oscillatory}. This unified approach to analyze and understand feedback-evolving games is called eco-evolutionary game theory. It denotes the coupled evolution of strategies and the environment. The cumulative feedback of decisions can subsequently alter environment-dependent incentives, thereby leading to new dynamical phenomena and new challenges for control.

Utilizing common resources always imposes a dilemma for community members. While cooperators restrain themselves and help maintain the proper state of resources, defectors demand more than their supposed share for a higher payoff. The nature of the feedback and the rate of ecological changes can relax or aggravate social dilemmas and promote persistent periodic oscillations of strategy abundance and environmental quality. To alleviate the tragedy of the commons, questions that inspire our present work include: how to manage and conserve public resources, which processes drive human cooperation, and how institutions, norms and other feedback mechanisms can be used to reinforce positive behaviors.

Strategies and the environment are coupled. In human societies the institutions structuring social interactions can be seen as part of the environment that co-evolve with strategic behaviors. Many problems need to be solved like pollution control and climate change mitigation, and the proper management of food webs, soil weathering, and earth system processes. The psychology of decision-making, species interaction and climate change action.

In this paper we show some simple examples where the original system has unstable equilibria under two distinct environments, and we can trap the system between the unstable equilibria to achieve a dynamic steady state for the system. We hope this can provide some intuition for policy making in social dilemmas. Controlled decision making is likely to be costly but will allow individuals to choose optimal behavior.

\section{Model}

Consider a two player asymmetrical game of normal form. Player 1 has the strategy set $S=\{s_1,s_2\}$, and player 2 has the strategy set $T=\{t_1,t_2\}$. When the pair of strategies $(s_i,t_j)$ is chosen, the payoff to player 1 is $a_{ij}=u_1(s_i,t_j)$ and the payoff to player 2 is $b_{ij}=u_2(s_i,t_j)$. Then the values of the payoff functions can be given by two matrices, respectively, for the two types of players 1 and 2,

$$A=\begin{bmatrix}a_{11} & a_{12} \\ a_{21} & a_{22}\end{bmatrix}, \quad B=\begin{bmatrix}
b_{11} & b_{12} \\
b_{21} & b_{22}
\end{bmatrix}.$$

We put these two matrices together in the following canonical way for a bimatrix game.

$$M=A\oplus B=\begin{bmatrix}
(a_{11},b_{11}) & (a_{12},b_{12}) \\
(a_{21},b_{21}) & (a_{22},b_{22})
\end{bmatrix}.$$

Now we consider a bimatrix game where the payoff bimatrix is not constant and depends on the environment.
$$M(t)=\begin{bmatrix}\left(a_{11}(t),b_{11}(t)\right) & \left(a_{12}(t),b_{12}(t)\right) \\
\left(a_{21}(t),b_{21}(t)\right) & \left(a_{22}(t),b_{22}(t)\right)\end{bmatrix}.$$
For simplicity, we start with the combination of two constant bimatrices,
\begin{align*}
    M(t)=\eta(t)M_I+(1-\eta(t))M_{II},
\end{align*}
where
\begin{align*}
    M_I=\begin{bmatrix}
    (a_{11},b_{11}) & (a_{12},b_{12}) \\
    (a_{21},b_{21}) & (a_{22},b_{22})
    \end{bmatrix},
    M_{II}=\begin{bmatrix}
    (a_{11}',b_{11}') & (a_{12}',b_{12}') \\
    (a_{21}',b_{21}') & (a_{22}',b_{22}')
    \end{bmatrix}
    ,\end{align*}
and $\eta(t)$ is the environment evolving indicator as the environment switches between two states I and II. There are many possibilities for the function $\eta(t)$. One of them is the following periodic function.

$$q(t)=\begin{cases}
0,\qquad kT\le t<kT+T/2, & k=0,1,2,\dots\\
1,\qquad kT+T/2\le t<(k+1)T, & k=0,1,2,\dots
\end{cases}$$

\begin{figure}[h]
    \centering
    \includegraphics[scale=0.37]{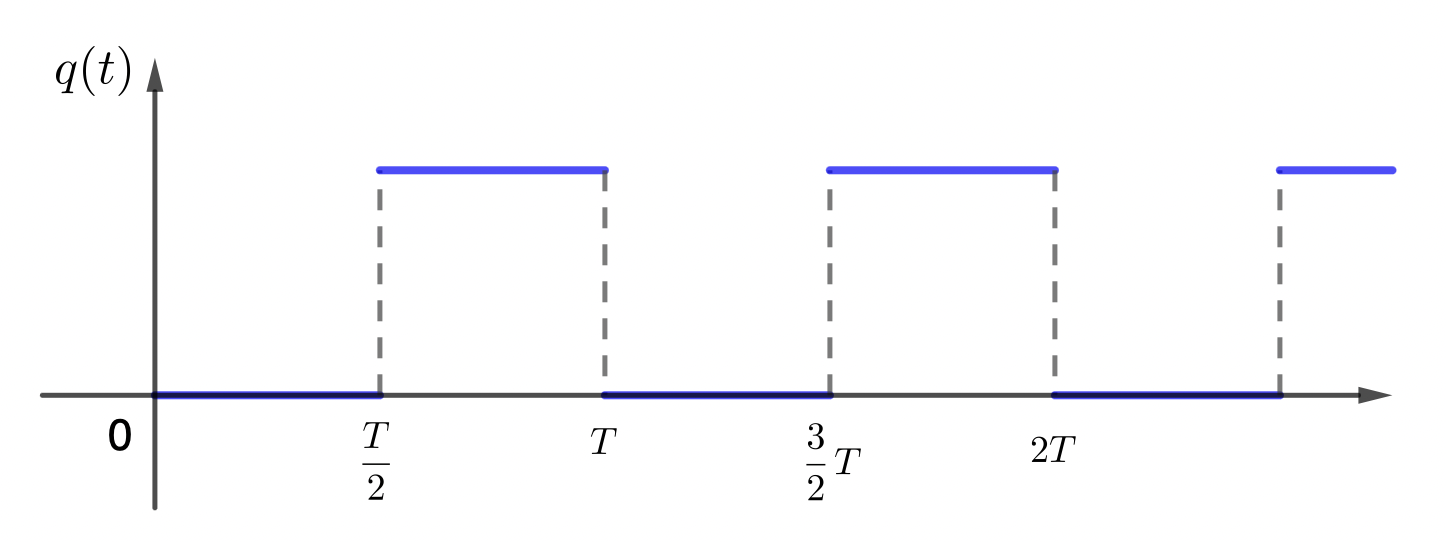}
    \caption{Periodic switching between two environments}
    \label{switch}
\end{figure}

Now we introduce mixed strategies in our game and assume $X=(x_1,x_2)$, where $x_1$ represents the probability of player 1 choosing strategy $s_1=l$ and $x_2$ represents the probability of player 1 choosing strategy $s_2=r$. We have $$x_1+x_2=1.$$ Similarly we have $Y=(y_1,y_2)$ where $y_1,y_2$ represent the probabilities of player 2 choosing strategies $t_1=L$ and $t_2=R$. Then
$$y_1+y_2=1.$$

We note that identical notations and equations can be used to describe eco-evolutionary dynamics in subpopulations consisting of type 1 and type 2 players. 

The payoff of player 1 when choosing strategy $s_1=l$ is
\begin{align*}
    \pi_l(X,Y)=a_{11}(t)y_1+a_{12}(t)y_2.
\end{align*}
The payoff of player 1 when choosing strategy $s_2=r$ is
\begin{align*}
    \pi_r(X,Y)=a_{21}(t)y_1+a_{22}(t)y_2.
\end{align*}
Let $E(\pi)=x_1\pi_l+x_2\pi_r$ be the expected payoff for player 1. The replicator equation for $x_1$ is
$$\dot x_1=x_1(\pi_l-E(\pi))=x_1(1-x_1)(\pi_l-\pi_r).$$
Plug in the values of $\pi_l,\pi_r$, we get
\begin{equation}
    \dot x_1=x_1(1-x_1)[(a_{11}(t)y_1+a_{12}(t)y_2)-(a_{21}(t)y_1+a_{22}(t)y_2)].
    \label{l}
\end{equation}
And the replicator equation for $x_2$ is
$$\dot x_2=-\dot x_1$$ since $x_1+x_2=1$.

Similarly we can write down the payoffs of strategies $t_1=L$ and $t_2=R$.
\begin{align*}
    &\pi_L(X,Y)=b_{11}(t)x_1+b_{21}(t)x_2,\\
    &\pi_R(X,Y)=b_{12}(t)x_1+b_{22}(t)x_2.
\end{align*}
The replicator equations for $y_1$ and $y_2$ are
\begin{align*}
    &\dot y_1=y_1(1-y_1)(\pi_L-\pi_R),\\
    &\dot y_2=-\dot y_1.
\end{align*}
Plug in the values of $\pi_L$ and $\pi_R$ we get
\begin{equation}
    \dot y_1=y_1(1-y_1)[(b_{11}(t)x_1+b_{21}(t)x_2)-(b_{12}(t)x_1+b_{22}(t)x_2)].
    \label{L}
\end{equation}

Equations \eqref{l},\eqref{L} combined with
\begin{equation}
    \eta(t)=q(t)
\end{equation}
governs the dynamics of our system.

\section{Simple One-Dimensional Model}
Before we try to solve the general 2-player bimatrix game $(A \neq B^T)$, let's start with the simplified one dimensional model obtained from assuming $A = B^T$.

In this case, we have a system governed by the replicator equation of the following form

\begin{equation}
    \dot x= x(1-x)[a(t)x-b(t)],\quad 0<x<1
    \label{1d-eq}
\end{equation}
where $a(t) = a_{11}(t) -a_{12}(t) - a_{21}(t) + a_{22}(t)$, and $b(t) = a_{22}(t) - a_{12}(t)$.

Let $a^*=b(t)/a(t)$. As one can see, this equation has three equilibrium points, $0,1, a^*$. For our purposes, we only consider cases where $a^*\in(0,1)$ by assuming the existence of an interior equilibrium,.

Based on the signs of $a(t)$ and $b(t)$, we have two possible situations for the phase line of the system.

\begin{figure}[h]
    \centering
    \includegraphics[scale=0.45]{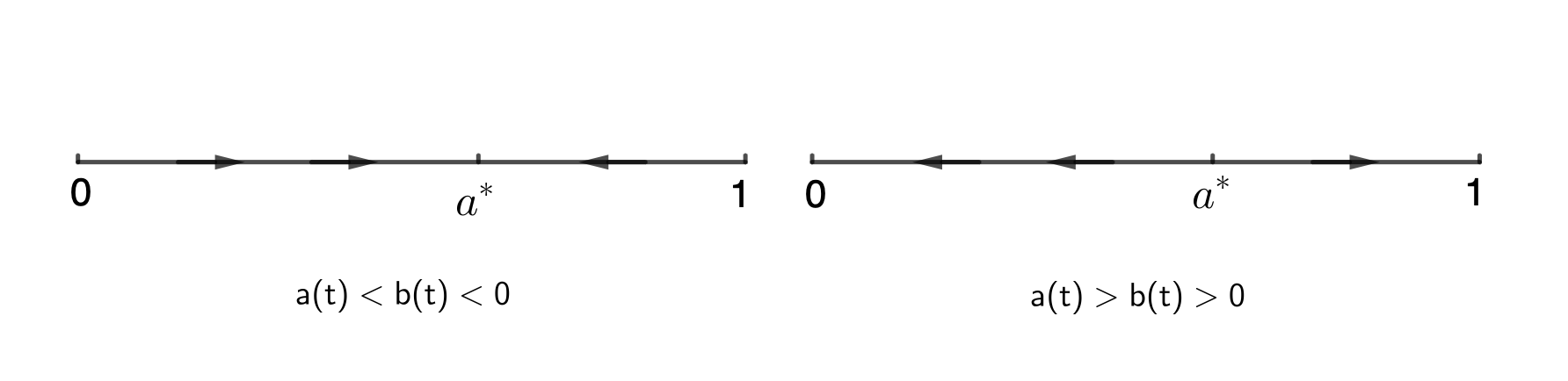}
    \caption{phase line of 1 dimensional replicator equation}
    \label{1-d}
\end{figure}

As Fig. \ref{1-d} shows, the dynamics of the system is quite clear if $a(t)$ and $b(t)$ are fixed constants. However if $a(t)$ and $b(t)$ are
changing over time, the dynamics becomes much more complicated.

Let's consider a simple case where $a(t)$ and $b(t)$ switch between two pairs of constants periodically. More specifically, we are interested in the case where the two separate states both have unstable fixed points, $a^*_1$ and $a^*_2$, yet if we choose the period carefully we can trap the system between $a^*_1$ and $a^*_2$, therefore achieving a dynamic equilibrium by switching between two different unstable systems.

Below we give two explicit examples with specific numbers for all the parameters to gain some intuition.

\subsection{Symmetric 1 dimensional example}
Consider two states governed by the following equations. By symmetry we mean that the term $a(t) = a_{11}(t) -a_{12}(t) - a_{21}(t) + a_{22}(t)$ is identical in both environments. 

\begin{equation}
    \dot x=x(1-x)(3x-1)
    \label{1d1}
\end{equation}
\begin{equation}
    \dot x=x(1-x)(3x-2)
    \label{1d2}
\end{equation}

Both equations are separable and we can solve them explicitly. In particular we care about the solution between the two unstable fixed points 1/3 and 2/3.

\begin{align*}
    \frac{dx}{x(1-x)(3x-1)}=(-\frac{1}{x}+\frac{1}{2}\frac{1}{1-x}+\frac{9}{2}\frac{1}{3x-1})dx=dt, \\
    -\ln x-\frac{1}{2}\ln(1-x)+\frac{3}{2}\ln(3x-1)=t+c, \\
    \frac{(3x-1)^3}{(1-x)x^2}=c_1e^{2t},\qquad  1/3<x<1.
\end{align*}

This is the solution to equation \eqref{1d1} and similarly we can solve \eqref{1d2}.

\begin{align*}
    \frac{dx}{x(1-x)(3x-2)}=(-\frac{1}{2}\frac{1}{x}+\frac{1}{1-x}+\frac{9}{2}\frac{1}{3x-2})dx=dt, \\
    -\frac{1}{2}\ln x-\ln(1-x)+\frac{3}{2}\ln(2-3x)=t+c, \\
    \frac{(2-3x)^3}{(1-x)^2x}=c_2e^{2t},\qquad 0<x<2/3.
\end{align*}

Assume at time $t=0$, $x$ starts at position $1/3+\epsilon,0<\epsilon<1/3$, and the system obeys equation \eqref{1d1}. Then as time goes on $x$ will move right toward $2/3$, if we manage to switch the governing equation to \eqref{1d2} before it reaches $2/3$, then $x$ will move back left.

Plug in the initial condition, we get

\begin{align*}
    c_1=\frac{(3\epsilon)^3}{(2/3-\epsilon)(1/3+\epsilon)^2}.
\end{align*}

Therefore $x$ will reach $2/3$ at time $t$ that satisfies

\begin{align*}
    c_1e^{2t}=&\frac{1}{1/3\times(2/3)^2}=27/4,\\
    t=&1/2[\ln(27/4)-\ln c_1] \\
    =&1/2[\ln(2/3-\epsilon)+2\ln(1/3+\epsilon)-3\ln\epsilon-2\ln2].
\end{align*}

Take the derivative of $t$ on $\epsilon$ we get

\begin{align*}
    t'=\frac{1}{2}[\frac{-1}{2/3-\epsilon}+\frac{2}{1/3+\epsilon}-\frac{3}{\epsilon}] <\frac{1}{2}[\frac{2}{1/3+\epsilon}-\frac{3}{\epsilon}]<0, \\0<\epsilon<1/3.
\end{align*}

Therefore $t$ decreases as $\epsilon$ increases between 0 and 1/3. The smallest time it takes for $x$ to go from its initial position to 2/3 happens when $\epsilon\to1/3$. 

In fact when $\epsilon=1/3$, we have

\begin{align*}
    t=1/2[\ln(1/3)+2\ln(2/3)+3\ln3-2\ln2]=0.
\end{align*}

And as $\epsilon\to0$, we can see $t\to\infty$. This means if we put $x$ initially at 1/3, it takes forever for $x$ to move over to 2/3. And if we put $x$ initially at 2/3, well, it's already there and if we don't want $x$ to move further right past 2/3, we have to immediately switch the system to satisfy equation \eqref{1d2}.

If we can control $\epsilon$, the upper bound of the time we want to switch over to equation \eqref{1d2} is given by the equation

\begin{equation}
    t_l=\frac{1}{2}[\ln(2/3-\epsilon)+2\ln(1/3+\epsilon)-3\ln\epsilon-2\ln2].
    \label{left1}
\end{equation}

Now let's look at the other end of our trapping zone, close to 2/3. Assume $x$ is positioned at $2/3-\delta,0<\delta<1/3$ at time $t=0$ and the system satisfies equation \eqref{1d2}.

Plug in the initial condition we get

\begin{align*}
    c_2=\frac{(3\delta)^3}{(1/3+\delta)^2(2/3-\delta)}.
\end{align*}

$x$ will move left and reach 1/3 at time $t$.

\begin{align*}
    c_2e^{2t}=\frac{1}{(2/3)^2\times1/3}=\frac{27}{4}.
\end{align*}

Unsurprisingly this is totally symmetric with the analysis we did before with $\epsilon$. Therefore we know here $t$ also decreases as $\delta$ increases between 0 and 1/3. In fact we can see the two differential equations are symmetric as well. If we can find the time it takes for $x$ to move from $1/3+\epsilon$ to $2/3-\epsilon$, $0<\epsilon<1/6$, we can trap $x$ in $(1/3+\epsilon,2/3-\epsilon)$ forever, just by switching between the two equations periodically. Now we compute critical value for the period $T$.

\begin{align*}
    c_1e^{2T}=\frac{(1-3\epsilon)^3}{(1/3+\epsilon)(2/3-\epsilon)^2}, \\
    e^{2T}=\left(\frac{1-3\epsilon}{3\epsilon}\right)^3\frac{1/3+\epsilon}{2/3-\epsilon},
\end{align*}
\begin{align*}
    T=\frac{1}{2}\{3[\ln(1-3\epsilon)-\ln(3\epsilon)]+\ln(1/3+\epsilon)-\ln(2/3-\epsilon)\}.
\end{align*}

Again we can take a derivative of $T$ on $\epsilon$,

\begin{align*}
    T'=\frac{1}{2}\{3[\frac{-3}{1-3\epsilon}-\frac{3}{3\epsilon}]+\frac{1}{1/3+\epsilon}-\frac{-1}{2/3-\epsilon}\} \\
    =\frac{1}{2}[(\frac{1}{1/3+\epsilon}-\frac{3}{\epsilon})+(\frac{1}{2/3-\epsilon}-\frac{3}{1/3-\epsilon})]<0,\\ 0<\epsilon<1/6.
\end{align*}

Again $T$ decreases as $\epsilon$ increases. When $\epsilon\to0$, $T\to\infty$. When $\epsilon=1/6,T=0$. For all values $\epsilon\in(0,1/6)$, we can achieve trapping between $(1/3+\epsilon,2/3-\epsilon)$ by setting the initial condition to be either $x(0)=1/3+\epsilon$ with equation \eqref{1d1} or $x(0)=2/3-\epsilon$ with equation \eqref{1d2}, and switching between the two equations with period

\begin{equation}
    T=\frac{1}{2}\{3[\ln(1-3\epsilon)-\ln(3\epsilon)]+\ln(1/3+\epsilon)-\ln(2/3-\epsilon)\}.
    \label{period1}
\end{equation}

\subsection{Asymmetric 1 dimensional example}

The above example is symmetric and therefore special. We'd like to explore the more general scenario of asymmetric cases, where $a(t) = a_{11}(t) -a_{12}(t) - a_{21}(t) + a_{22}(t)$ is of different values under each environment.

Let's replace equation \eqref{1d1} with the following equation.

\begin{equation}
    \dot x=x(1-x)(4x-1).
    \label{1d3}
\end{equation}

The unstable equilibrium of equation \eqref{1d3} happens at 1/4. We want to study the behavior of the solution between 1/4 and 2/3.

\begin{align*}
    \frac{dx}{x(1-x)(4x-1)}=(-\frac{1}{x}+\frac{1}{3}\frac{1}{1-x}+\frac{16}{3}\frac{1}{4x-1})dx=dt, \\
    -\ln x-\frac{1}{3}\ln(1-x)+\frac{4}{3}\ln(4x-1)=t+c, \\
    \frac{(4x-1)^4}{x^3(1-x)}=c_3e^{3t}, \qquad 1/4<x<1.
\end{align*}

Again assume we start the system at $x(0)=1/4+\epsilon, 0<\epsilon<5/12$ with equation \eqref{1d3}, then we get

\begin{align*}
    c_3=\frac{(4\epsilon)^4}{(1/4+\epsilon)^3(3/4-\epsilon)}.
\end{align*}

And $x$ arrives at 2/3 at time $t$, then

\begin{align*}
    c_3e^{3t}=&\frac{(5/3)^4}{(2/3)^3(1/3)}=\frac{5^4}{2^3},\\
    t=&\frac{1}{3}[4\ln5-3\ln2-\ln c_3]\\
    =&\frac{1}{3}\{4[\ln 5-\ln(4\epsilon)]+3[\ln(1/4+\epsilon)-\ln2]+\ln(3/4-\epsilon)\}.
\end{align*}

Again we can take a derivative,

\begin{align*}
    t'=\frac{1}{3}(-\frac{4}{\epsilon}+\frac{3}{1/4+\epsilon}-\frac{1}{3/4-\epsilon})<0,\\ 0<\epsilon<5/12.
\end{align*}

As $\epsilon\to0$, $t\to\infty$ since 1/4 is an equilibrium point. And  $t=0$ if $\epsilon=5/12$. The time it takes for $x$ to move from $1/4+\epsilon$ to $2/3-\delta, 0<\delta<5/12-\epsilon<5/12$ is

\begin{equation}
    t_l=\frac{1}{3}\{4[\ln(5/3-4\delta)-\ln(4\epsilon)]+3[\ln(1/4+\epsilon)-\ln(2/3-\delta)]+\ln(3/4-\epsilon)-\ln(1/3+\delta)\}.
    \label{left2}
\end{equation}

Now if we switch the system to equation \eqref{1d2}, the time it takes for $x$ to move from $2/3-\delta$ back to $1/4+\epsilon$ is

\begin{equation}
    t_r=\frac{1}{2}\{3[\ln(5/4-3\epsilon)-\ln(3\delta)]+2[\ln(1/3+\delta)-\ln(3/4-\epsilon)]+\ln(2/3-\delta)-\ln(1/4+\epsilon)\}.
    \label{right2}
\end{equation}

There is no obvious reason to think $t_l$ and $t_r$ are equal. In fact, if we take $\epsilon=1/12,\delta=1/6$, we can compute $t_l=5/3\ln2\neq1/2\ln(27/4)=t_r$.

So under this asymmetric system, we could start with equation \eqref{1d3} at $x(0)=1/4+\epsilon$, switch to equation \eqref{1d2} after time $t_l$, then switch back to equation \eqref{1d3} after time $t_r$ and repeat. Or we could start the system with equation \eqref{1d2} at $x(0)=2/3-\delta$, switch to equation \eqref{1d3} after time $t_r$, then switch back to equation \eqref{1d2} after time $t_l$ and repeat. In either case, we manage to trap our system between $1/4+\epsilon$ and $2/3-\delta$.

With this asymmetric example, the time periods we stay on each state is different, but we can still achieve trapping by repeating a fixed control pattern of staying on one system for a fixed period and then the other on a different fixed period.

\subsection{General asymmetric one-dimensional model}

In general, let's say we have the system switch between the following two equations

\begin{equation}
    \dot x=x(1-x)(a_1x-b_1), a_1>b_1>0.
    \label{1dgl}
\end{equation}

\begin{equation}
    \dot x=x(1-x)(a_2x-b_2), a_2>b_2>0.
    \label{1dgr}
\end{equation}

Without loss of generality, we assume $a_1^*=b_1/a_1<a_2^*=b_2/a_2$. Again we want to trap the system between $(a_1^*,a_2^*)\subsetneq(0,1)$.

Assume we start the system on equation \eqref{1dgl} at $x(0)=a_1^*+\epsilon$, and we switch the system to equation \eqref{1dgr} when $x$ reaches $a_2^*-\delta$, then we switch back to equation \eqref{1dgl} after $x$ reaches $a_1^*+\epsilon$ again and repeat this process. Here $\epsilon,\delta>0,\epsilon+\delta<a_2^*-a_1^*$. To find the period the system stays on equation \eqref{1dgl}, we need to solve the equation.

\begin{align*}
    \frac{dx}{x(1-x)(a_1x-b_1)}=[-\frac{1}{b_1}\frac{1}{x}+\frac{1}{a_1-b_1}\frac{1}{1-x}+\frac{a_1^2}{b_1(a_1-b_1)}\frac{1}{a_1x-b_1}]dx=dt, \\
    -\frac{1}{b_1}\ln x-\frac{1}{a_1-b_1}\ln(1-x)+\frac{a_1}{b_1(a_1-b_1)}\ln(a_1x-b_1)=t+c, \\
    \frac{(a_1x-b_1)^{a_1}}{x^{a_1-b_1}(1-x)^{b_1}}=c_le^{b_1(a_1-b_1)t}, \qquad a_1^*<x<1.
\end{align*}

Plug in the initial condition, we get

\begin{align*}
    c_l=\frac{(a_1\epsilon)^{a_1}}{(a_1^*+\epsilon)^{a_1-b_1}(1-a_1^*-\epsilon)^{b_1}}.
\end{align*}

The time it takes for $x$ to move to $a_2^*-\delta$ is

\begin{multline}
    t_l=\frac{1}{b_1(a_1-b_1)}\{a_1[\ln(a_1a_2^*-a_1\delta-b_1)-\ln(a_1\epsilon)]+ \\
    (a_1-b_1)[\ln(a_1^*+\epsilon)-\ln(a_2^*-\delta)]+ \\
    b_1[\ln(1-a_1^*-\epsilon)-\ln(1-a_2^*+\delta)]\}.
    \label{left}
\end{multline}

One can check formulas \eqref{left1}, \eqref{period1}, and \eqref{left2} from the symmetric and asymmetric examples are special cases of \eqref{left}.

Similarly we can solve \eqref{1dgr}.

\begin{align*}
    \frac{dx}{x(1-x)(a_2x-b_2)}=[-\frac{1}{b_2}\frac{1}{x}+\frac{1}{a_2-b_2}\frac{1}{1-x}+\frac{a_2^2}{b_2(a_2-b_2)}\frac{1}{a_2x-b_2}]dx=dt, \\
    -\frac{1}{b_2}\ln x-\frac{1}{a_2-b_2}\ln(1-x)+\frac{a_2}{b_2(a_2-b_2)}\ln(b_2-a_2x)=t+c, \\
    \frac{(b_2-a_2x)^{a_2}}{x^{a_2-b_2}(1-x)^{b_2}}=c_re^{b_2(a_2-b_2)t}, \qquad 0<x<a_2^*.
\end{align*}

Plug in the initial conditon $x(0)=a_2^*-\delta$, we get

\begin{align*}
    c_r=\frac{(a_2\delta)^{a_2}}{(a_2^*-\delta)^{a_2-b_2}(1-a_2^*+\delta)^{b_2}}.
\end{align*}

The time it takes $x$ to move from $a_2^*-\delta$ to $a_1^*+\epsilon$ is

\begin{multline}
    t_r=\frac{1}{b_2(a_2-b_2)}\{a_2[\ln(b_2-a_2a_1^*-a_2\epsilon)-\ln(a_2\delta)]+ \\
    (a_2-b_2)[\ln(a_2^*-\delta)-\ln(a_1^*+\epsilon)]+ \\
    +b_2[\ln(1-a_2^*+\delta)-\ln(1-a_1^*-\epsilon)]\}.
    \label{right}
\end{multline}

Again we can check formulas \eqref{period1} and \eqref{right2} are special cases of \eqref{right}.

We can trap our system between $(a_1^*,a_2^*)$ by starting from equation \eqref{1dgl} at $x(0)=a_1^*+\epsilon$, switch the equation to \eqref{1dgr} after time $t_l$, then switch back to \eqref{1dgl} after time  $t_r$, and repeat this process. Or we can start from equation \eqref{1dgr} with $x(0)=a_2^*-\delta$, switch the equation after time $t_r$ to \eqref{1dgl}, and switch again back to \eqref{1dgr} after time $t_l$ and repeat the process. Here $t_r$ and $t_l$ are given by formulas \eqref{left}, \eqref{right}.

\subsection{Continuous variation}

Both the symmetric and asymmetric examples above switch between two different equations. What if the equation parameters vary continuously?

Still assuming $a(t)>b(t)>0$, let $a^*=b(t)/a(t)$ be a periodic function that ranges between $p_0$ and $p_1$, $0<p_0<p_1<1$. Again we can analyze the behavior of the system with figure \ref{1-d}. One way to ensure $x(t)$ is trapped between $p_0$ and $p_1$, at any moment in time, is if $a^*$ is increasing, we want $x(t)>a^*$ so that $x(t)$ doesn't come to the left of $a^*$ and converge to 0. Similarly if $a^*$ is decreasing, we want $x(t)<a^*$ so that $x(t)$ doesn't converge to 1.

To summarize, we need

\begin{equation}
    [x(t)-a^*]\dot a^*\ge0,p_0<x(t)<p_1,\qquad\forall t.
    \label{continuous}
\end{equation}

for the system to be trapped in between $p_0$ and $p_1$. Plug $a^*=b(t)/a(t)$ into \eqref{continuous}, we get

\begin{equation}
    [a(t)x(t)-b(t)][b'(t)a(t)-a'(t)b(t)]\ge0,p_0<x(t)<p_1,\qquad\forall t.
\end{equation}

\section{Dynamics of the Mixed Strategy Bimatrix Game}
Considering $x_1+x_2=1,y_1+y_2=1$, our system of mixed strategy bimatrix game has the following replicator equations
$$
\begin{cases}
\dot x_1=x_1(1-x_1)\{[a_{11}(t)+a_{22}(t)-a_{12}(t)-a_{21}(t)]y_1-[a_{22}(t)-a_{12}(t)]\},\\
\dot y_1=y_1(1-y_1)\{[b_{11}(t)+b_{22}(t)-b_{12}(t)-b_{21}(t)]x_1-[b_{22}(t)-b_{21}(t)]\}.
\end{cases}
$$
Here 
\begin{align*}
M(t)&=\begin{bmatrix}\left(a_{11}(t),b_{11}(t)\right) & \left(a_{12}(t),b_{12}(t)\right) \\
\left(a_{21}(t),b_{21}(t)\right) & \left(a_{22}(t),b_{22}(t)\right)\end{bmatrix}\\
&=\begin{cases}
M_{II}, t\in[kT,kT+T/2),\\
M_I, t\in[kT+T/2,(k+1)T),
\end{cases} k=0,1,2,\dots\\
M_I&=\begin{bmatrix}
    (a_{11},b_{11}) & (a_{12},b_{12}) \\
    (a_{21},b_{21}) & (a_{22},b_{22})
    \end{bmatrix},
    M_{II}=\begin{bmatrix}
    (a_{11}',b_{11}') & (a_{12}',b_{12}') \\
    (a_{21}',b_{21}') & (a_{22}',b_{22}')
    \end{bmatrix}.
\end{align*}

Although we started with a 2 player bimatrix game with mixed strategies, where $x_1$ and $y_1$ represent the probabilities of each player choosing the corresponding strategies, this model also works for two groups of players. There $x_1$ and $y_1$ represent the proportions or percentages of each player group that adopt the corresponding strategies.

Let $$a^*=\frac{a_{22}(t)-a_{12}(t)}{a_{11}(t)+a_{22}(t)-a_{12}(t)-a_{21}(t)},$$
$$b^*=\frac{b_{22}(t)-b_{21}(t)}{b_{11}(t)+b_{22}(t)-b_{12}(t)-b_{21}(t)}.$$ Then $(b^*,a^*)$ is a fixed point of our system, in addition to the points $(0,0)$, $(1,1)$, $(0,1)$, and $(1,0)$. Since $x_1$ and $y_1$ are probabilities or proportions, we only consider their values in $[0,1]$. The reasonable range for $a^*$ and $b^*$ is $(0,1)$. For example, for times $kT\le t<kT+T/2, k=0,1,2,\dots$, if $A'=\begin{bmatrix}
a_{11}' & a_{12}' \\
a_{21}' & a_{22}'
\end{bmatrix}$ and $B'=\begin{bmatrix}
b_{11}' & b_{12}' \\
b_{21}' & b_{22}'
\end{bmatrix}$ are diagonal matrices and their diagonal entries are positive,  we have
$$a^*=\frac{a_{22}'}{a_{11}'+a_{22}'}\in(0,1),\quad b^*=\frac{b_{22}'}{b_{11}'+b_{22}'}\in(0,1).$$

\begin{figure}[h]
    \centering
    \includegraphics[scale=0.4]{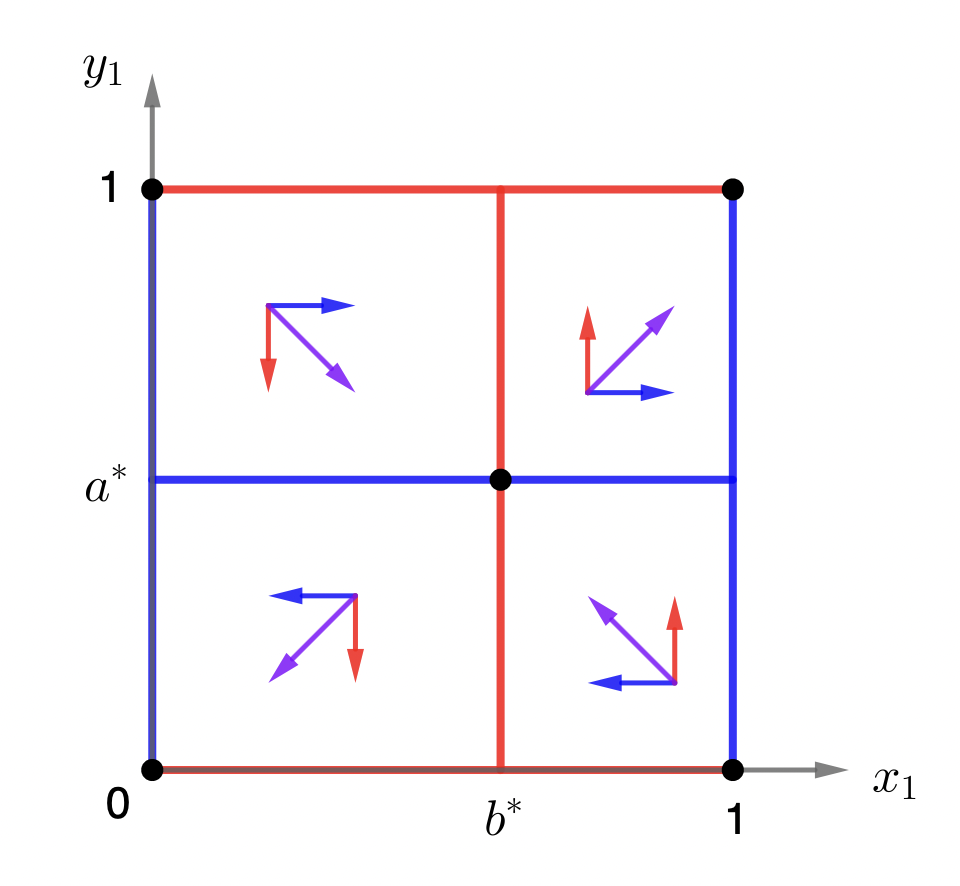}
    \caption{Phase plot for positive diagonal matrices}
    \label{diagonal}
\end{figure}

To analyze the dynamics of the system, we plot the nullclines and phase plot below (Fig. \ref{diagonal}). Here the blue nullclines are for $x_1$ and the red ones are for $y_1$. The intersection of a blue line and a red line is a fixed point, or equilibrium point. There are 5 equilibirum points in total. The blue arrows indicate whether $x_1$ is increasing or decreasing in the region, and the red arrows indicate the behavior of $y_1$ value over time. The purple arrow indicates qualitatively the direction of movements for points in that region. According to this qualitative phase portrait, we can see the equilibria $(0,0),(1,1)$ are locally stable, and equilibria $(0,1),(1,0)$ are unstable. The coexistence equilibrium $(b^*, a^*)$ is unstable. To be more specific, it is a saddle point.

If the matrices $A=\begin{bmatrix} a_{11} & a_{12} \\ a_{21} & a_{22} \end{bmatrix}$ and $B=\begin{bmatrix}
b_{11} & b_{12} \\
b_{21} & b_{22}
\end{bmatrix}$ are also positive diagonal matrices, then for times $kT+T/2\le t<(k+1)T, k=0,1,2,\dots$, we have the same qualitative phase portrait as Fig. \ref{diagonal}.

\begin{figure}[h]
    \centering
    \includegraphics[scale=0.4]{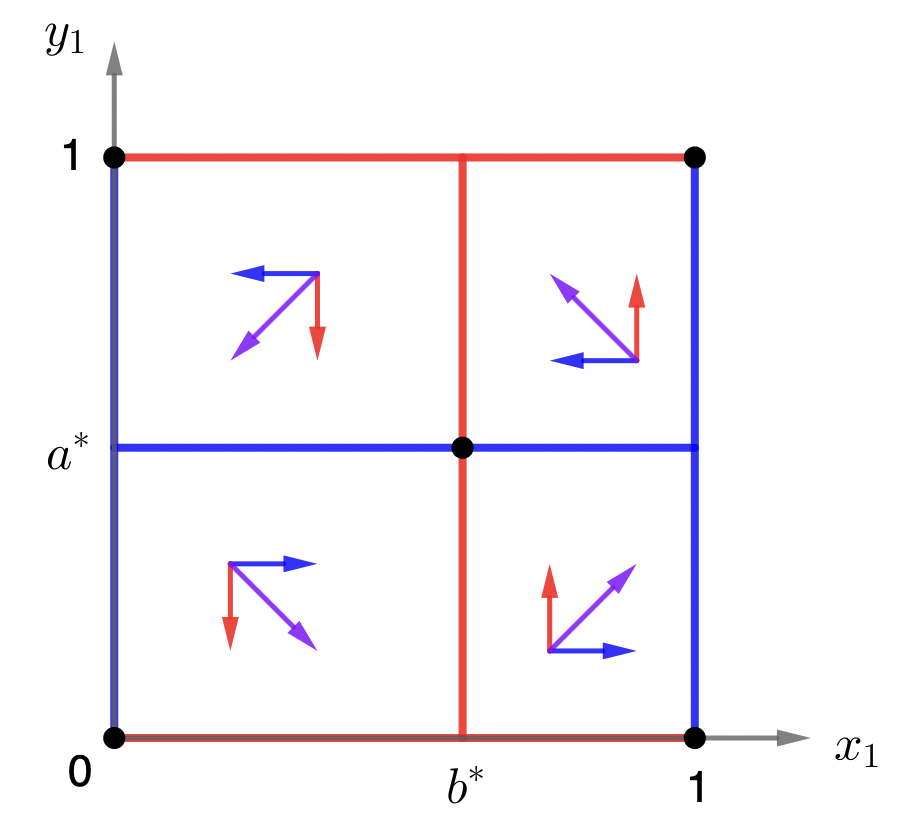}
    \caption{Oscillation around coexistence}
    \label{osc}
\end{figure}

If we want the probabilities/proportions to oscillate around the coexistence equilibrium $(b^*,a^*)$, the qualitative phase portrait should look like Fig. \ref{osc}.

For such a phase portrait, we need $0<a^*,b^*<1$ and also 
$$a_{11}(t)+a_{22}(t)-a_{12}(t)-a_{21}(t)<0,$$
$$b_{11}(t)+b_{22}(t)-b_{12}(t)-b_{21}(t)>0.$$

Solve all four inequalities, we get
$$a_{11}(t)<a_{21}(t),a_{22}(t)<a_{12}(t);$$
$$b_{11}(t)>b_{12}(t),b_{22}(t)>b_{21}(t).$$

These are the conditions for an oscillatory system around the coexistence equilibrium. In this case the coexistence equilibrium is a center.

\section{Trapping Regions}

Assume we have two different saddle points under environment matrices $M_I$ and $M_{II}$. We wish to do what we did with the one dimensional examples, achieve trapping between the two unstable fixed points by switching between the two environments.

Let's draw a qualitative phase portrait. Assume $(0,0)$ and $(1,1)$ are always stable, we have $E_1$ as the saddle point for environment $M_I$ and $E_2$ as the saddle point for environment $M_{II}$. 

\begin{figure}[h]
    \centering
    \includegraphics[scale=0.4]{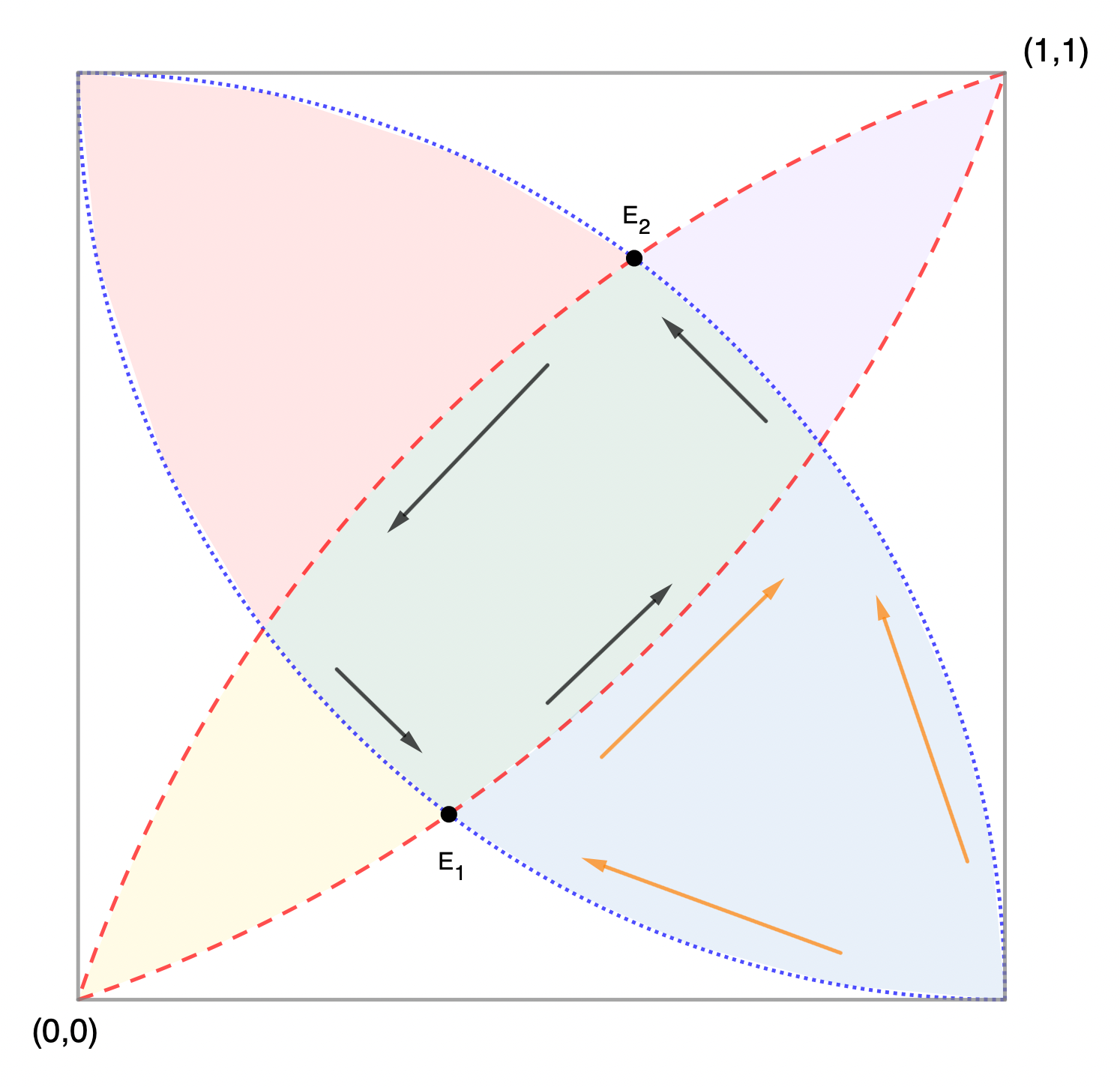}
    \caption{Trapping Region}
    \label{trapping2d}
\end{figure}

Here the blue dense dotted lines represent the stable manifold for the saddle points and the red sparse dotted lines represent the unstable manifolds. Then as long as the blue lines and red lines form a region bounded by alternating colored edges, we have a trapping region, shaded in green in Fig. \ref{trapping2d}. This is because, qualitatively, under environment I, the movements of points inside the region are marked by the two lower arrows, and when we switch to environment II, the movements of points inside the region are marked by the two upper arrows. Together they form a closed circular motion. As long as we make sure to switch the environment before the points moves out the boundary of the region, they should be trapped in the region forever.

Now let's look at the blue shaded region instead. Unlike the green region, along the boundaries of the blue region, the three orange arrows don't form a consecutive circular path. Indeed, when under environment I, the movements of points inside the blue region should follow the two orange arrows to the left, and if we don't ever switch the environment, the points will escape the blue region from the right. Now if we switch the environment before the points escape, they will now follow the orange arrow on the right. If we stay on environment II for long enough the points will move on into the green region and escape the blue region. If we manage to switch back to environment I before points move into the green region, well, they will keep moving right until they escape the blue region on the right.

This analysis can be applied to the purple, pink, and yellow regions as well. We can have a trapping region if and only if the arrows along the boundary form a closed circular path.

\section{Linearization of the system}

Although we don't have explicit solutions for our nonlinear system, we can still linearize it near the fixed point $(b^*,a^*)$.

\begin{equation}
    \begin{bmatrix}
    x_1 \\ y_1
    \end{bmatrix}'=\begin{bmatrix}
    0 & \alpha \\
    \beta & 0
    \end{bmatrix}(\begin{bmatrix}
    x_1 \\ y_1
    \end{bmatrix}-\begin{bmatrix}
    b^* \\ a^*
    \end{bmatrix}),
\end{equation}

here $$\alpha=b^*(1-b^*)[a_{11}(t)+a_{22}(t)-a_{12}(t)-a_{21}(t)],$$
$$\beta=a^*(1-a^*)[b_{11}(t)+b_{22}(t)-b_{12}(t)-b{21}(t)].$$

Still assuming we have a system that switches between environments \\
$M_I=A\oplus B$ and $M_{II}=A'\oplus B'$. Here $A=\begin{bmatrix}
a_{11} & a_{12} \\
a_{21} & a_{22}
\end{bmatrix},B=\begin{bmatrix}
b_{11} & b_{12} \\
b_{21} & b_{22}
\end{bmatrix},\\A'=\begin{bmatrix}
a_{11}' & a_{12}' \\
a_{21}' & a_{22}'
\end{bmatrix},B'=\begin{bmatrix}
b_{11}' & b_{12}' \\
b_{21}' & b_{22}'
\end{bmatrix}'$ are all constant matrices. Now assume 

\begin{align*}
    & a^*=\frac{a_{22}-a_{12}}{a_{11}+a_{22}-a_{12}-a_{21}}, \\
    & b^*=\frac{b_{22}-b_{21}}{b_{11}+b_{22}-b_{12}-b_{21}}, \\
    & c^*=\frac{a_{22}'-a_{12}'}{a_{11}'+a_{22}'-a_{12}'-a_{21}'},\\
    & d^*=\frac{b_{22}'-b_{21}'}{b_{11}'+b_{22}'-b_{12}'-b_{21}'};\\
    & \alpha=b^*(1-b^*)(a_{11}+a_{22}-a_{12}-a_{21}), \\
    & \beta =a^*(1-a^*)(b_{11}+b_{22}-b_{12}-b_{21}), \\
    & \alpha'=d^*(1-d^*)(a_{11}'+a_{22}'-a_{12}'-a_{21}'),\\
    & \beta'=c^*(1-c^*)(b_{11}'+b_{22}'-b_{12}'-b_{21}').
\end{align*}

Under environment $M_I$, the linear system of ODEs has solutions

$$\begin{bmatrix}
x_1 \\ y_1
\end{bmatrix}=\begin{bmatrix}
b^* \\ a^*
\end{bmatrix} + c_1e^{\sqrt{\alpha\beta}t}\begin{bmatrix}
\sqrt{\alpha} \\ \sqrt{\beta}
\end{bmatrix} + c_2e^{-\sqrt{\alpha\beta}t}\begin{bmatrix}
\sqrt{\alpha} \\ -\sqrt{\beta}
\end{bmatrix}.$$

And under environment $M_{II}$, this linear system of ODEs has solutions

$$\begin{bmatrix}
x_1 \\ y_1
\end{bmatrix}=\begin{bmatrix}
d^* \\ c^*
\end{bmatrix} + c_1e^{\sqrt{\alpha'\beta'}t}\begin{bmatrix}
\sqrt{\alpha'} \\ \sqrt{\beta'}
\end{bmatrix} + c_2e^{-\sqrt{\alpha'\beta'}t}\begin{bmatrix}
\sqrt{\alpha'} \\ -\sqrt{\beta'}
\end{bmatrix}.$$

\begin{figure}[h]
    \centering
    \includegraphics[scale=0.5]{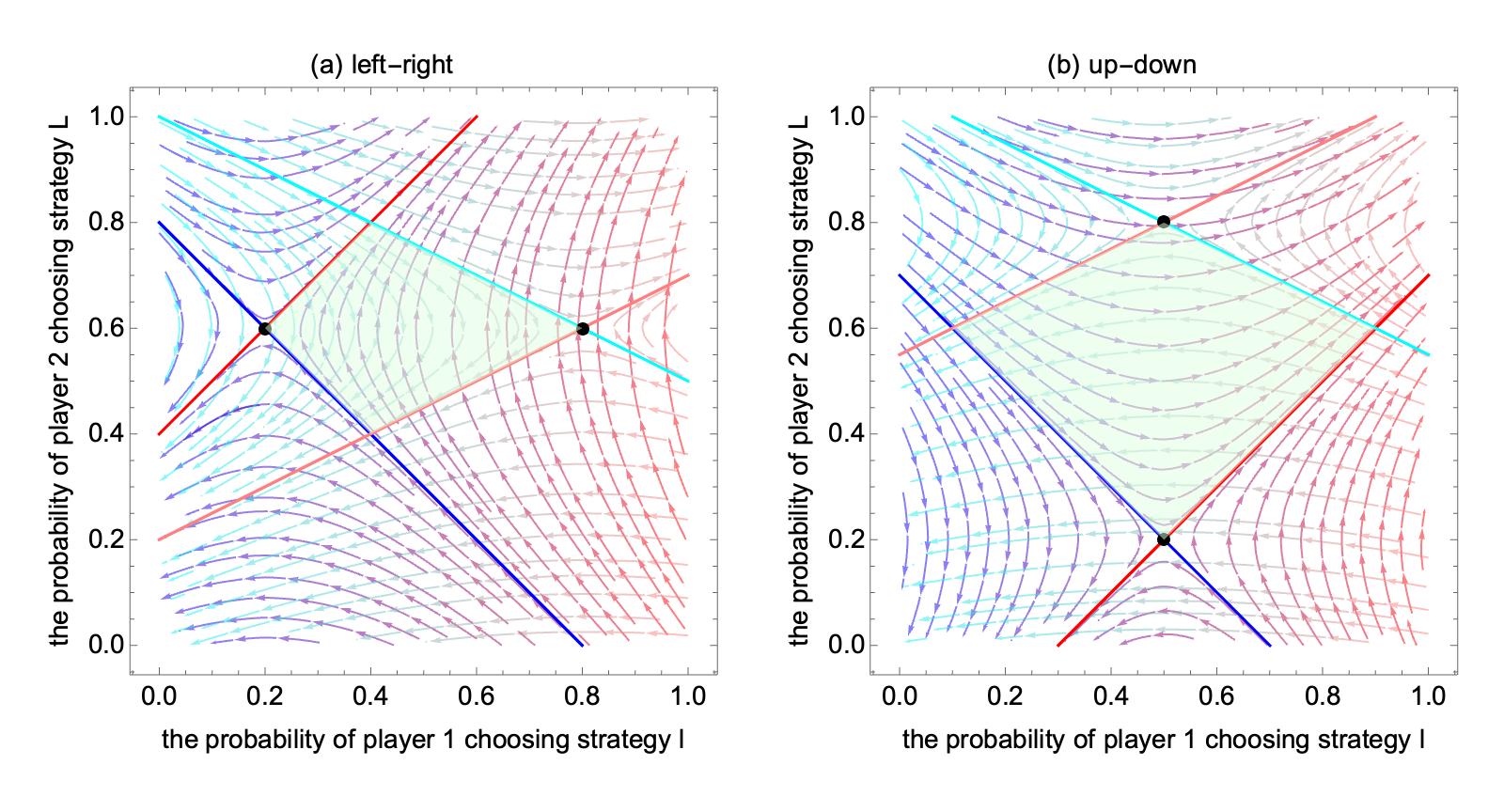}
    \caption{Linear Trapping: no shared manifolds}
    \label{left-right}
\end{figure}

The fixed points $E_1=(b^*,a^*),E_2=(d^*,c^*)$ are saddle points. Notice that their two eigenvector directions are symmetric by the $y$-axis. If the two fixed points are positioned as in Fig. \ref{left-right}(a), then the green shaded region is a trapping region. Here the blue and light blue straight lines represent the contracting directions, or stable manifolds and the red and pink lines represent expanding directions, or unstable manifolds.

To have this left-right configuration, we need

\begin{align}
    \left|\frac{a^*-c^*}{b^*-d^*}\right|<\min\{\sqrt{\frac{\beta}{\alpha}},\sqrt{\frac{\beta'}{\alpha'}}\}.
    \label{LR}
\end{align}

The contracting and expanding directions of $E_1$ divides our phase space into 4 parts. We should have 4 configurations with a convex quadrilateral trapping region. Only 2 of these are distinct considering symmetry between $E_1$ and $E_2$, we can always rename or relabel them. The other distinct configuration, the up-down configuration in Fig. \ref{left-right} (b) requires
\begin{align}
    \left|\frac{a^*-c^*}{b^*-d^*}\right|>\max\{\sqrt{\frac{\beta}{\alpha}},\sqrt{\frac{\beta'}{\alpha'}}\}.
    \label{UD}
\end{align}

\begin{figure}[H]
    \centering
    \includegraphics[scale=0.5]{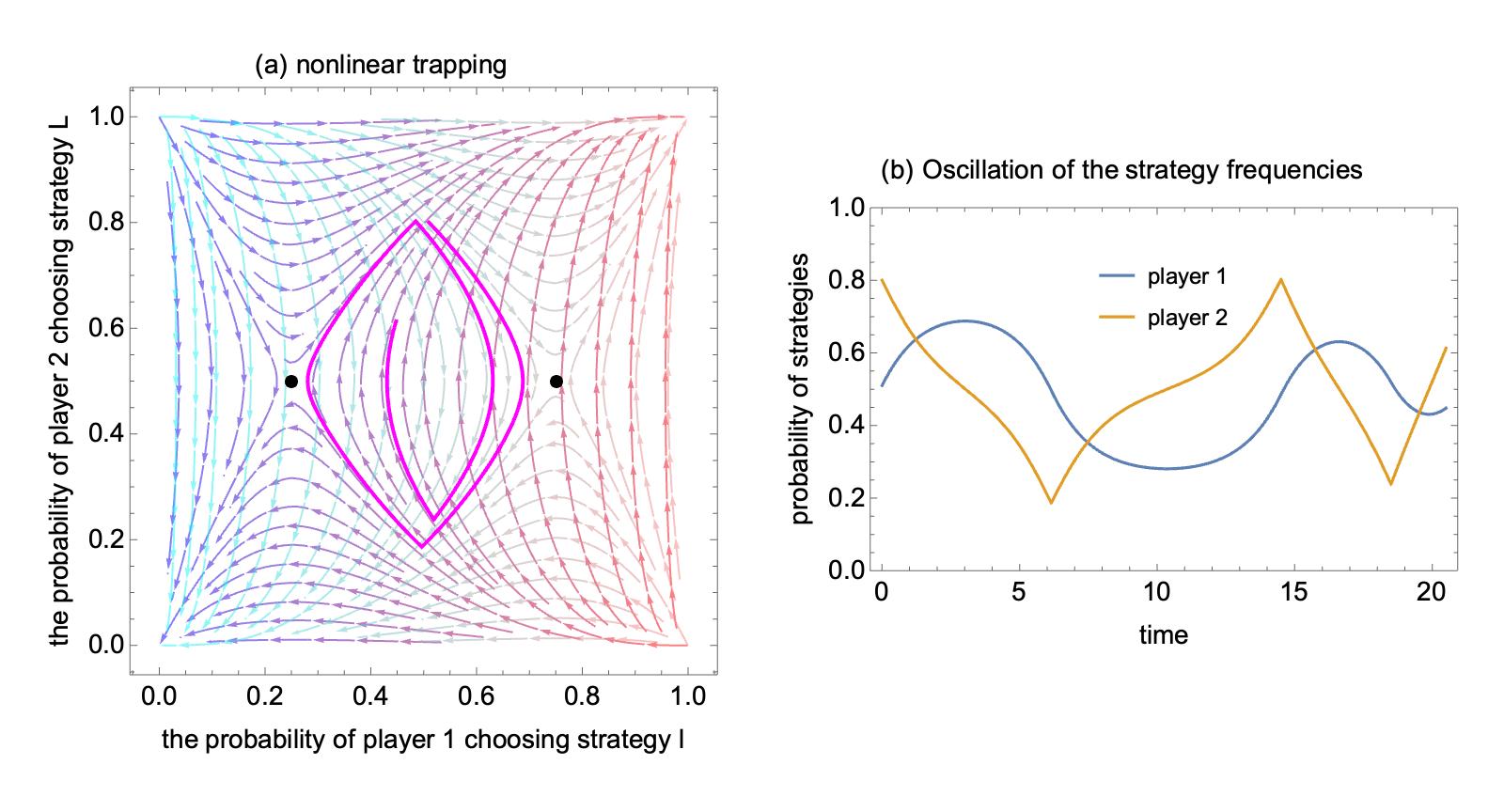}
    \caption{Nonlinear Trapping}
    \label{nonlinear}
\end{figure}

To show that this linear approximation actually is close enough to the nonlinear system, we study a specific nonlinear system here as well. The trapping region in Fig. \ref{nonlinear}(a) is not a convex quadrilateral any more, we plot a specific trapping trajectory in magenta and the corresponding frequencies of strategies in Fig.\ref{nonlinear}(b). The equations for this nonlinear system are

\begin{equation}
    \begin{cases}
    x'=x(1-x)(2y-1),\\
    y'=y(1-y)(4x-3).
    \end{cases}
    \label{non1}
\end{equation}
\begin{equation}
    \begin{cases}
    x'=x(1-x)(2y-1),\\
    y'=y(1-y)(4x-1).
    \end{cases}
    \label{non2}
\end{equation}

At time $t_0=0$, we start with \eqref{non1}, and initial conditions $x_0=0.51,y_0=0.8$. We switch over to \eqref{non2} at $t_1=6.15$. Then after $t_2=8.35$ we switch again, and the trajectory in magenta happens when we switch again after $t_3=4$ and $t_4=2$. The corresponding solutions $x(t)$ and $y(t)$ are plotted here as well (Fig. \ref{nonlinear}). We can see the curve for $x(t)$ is rather smooth while the curve for $y(t)$ has sharp tips. This is because the equation for $x(t)$ actually stays the same between the two environments. 

This example is rather simple but it is nonlinear, and it corresponds very well to the left-right linear trapping region (Fig. \ref{left-right} (a)). Given this resemblance, we have a good reason to think the linear trapping regions, although not exactly what we will see in realistic nonlinear systems, still give us reasonable approximation in nonlinear trapping.

In addition to the left-right and up-down configurations in Fig. \ref{left-right}, if $E_2$ lies on the stable or unstable manifold of $E_1$, or the other way around, we will also have a trapping region. These cases, co require

\begin{align}
    \left|\frac{a^*-c^*}{b^*-d^*}\right|=\sqrt{\frac{\beta}{\alpha}}\text{ or }\sqrt{\frac{\beta'}{\alpha'}}.
    \label{SUM}
\end{align}

\begin{figure}[h]
    \centering
    \includegraphics[scale=0.5]{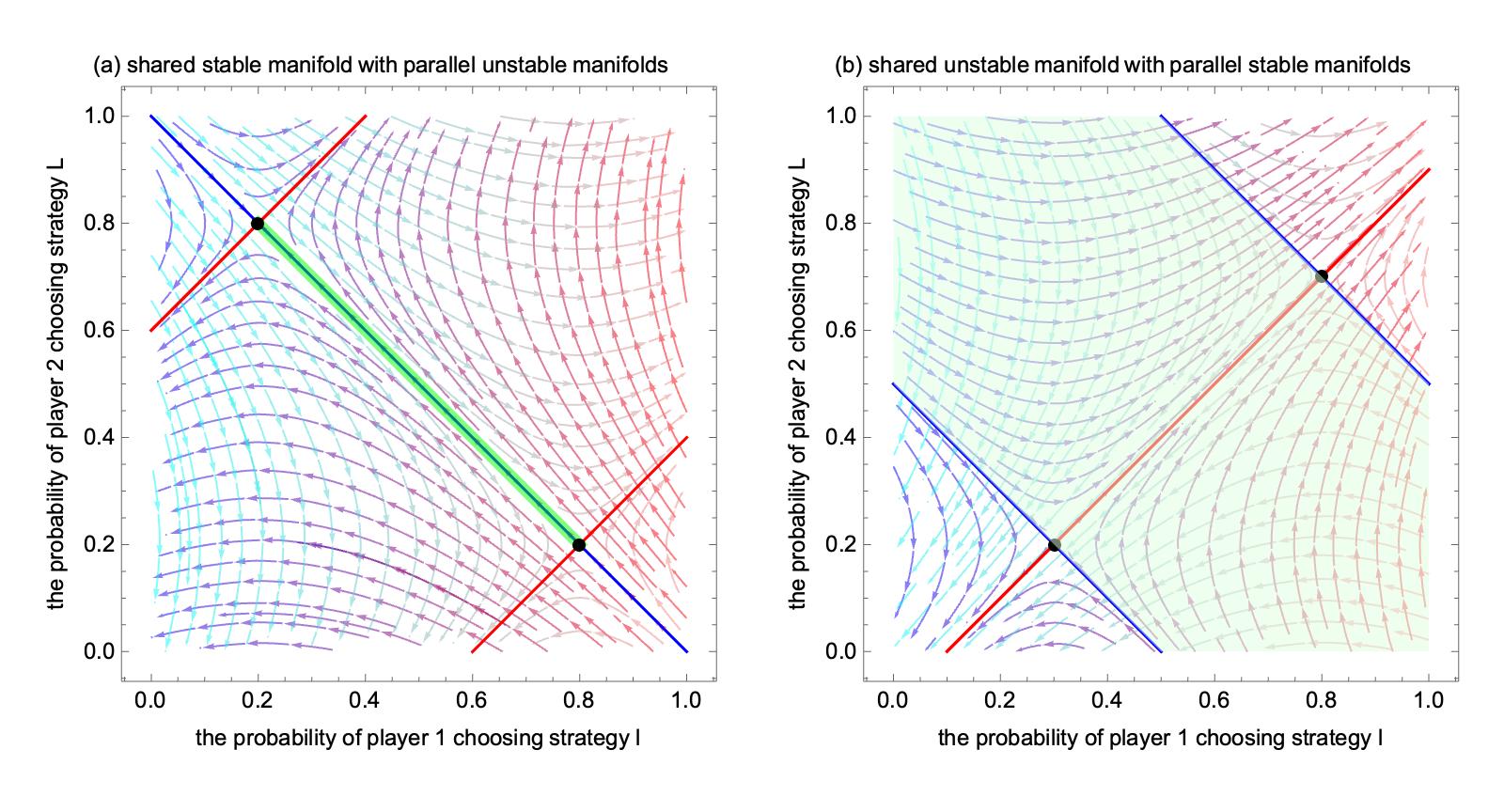}
    \caption{Linear Trapping: shared manifold with parallel other manifolds}
    \label{shared}
\end{figure}

In Fig. \ref{shared}, we have the very special case where one fixed point is on the stable, or unstable manifold of the other fixed point, and their other manifolds are parallel to each other. When one fixed point is on the stable manifold of the other, the trapping region degenerates to a line segment between the two fixed points. While when one fixed point is on the unstable manifold of the other, the trapping region becomes much larger and takes up a big part of our phase space.

\begin{figure}[H]
    \centering
    \includegraphics[scale=0.5]{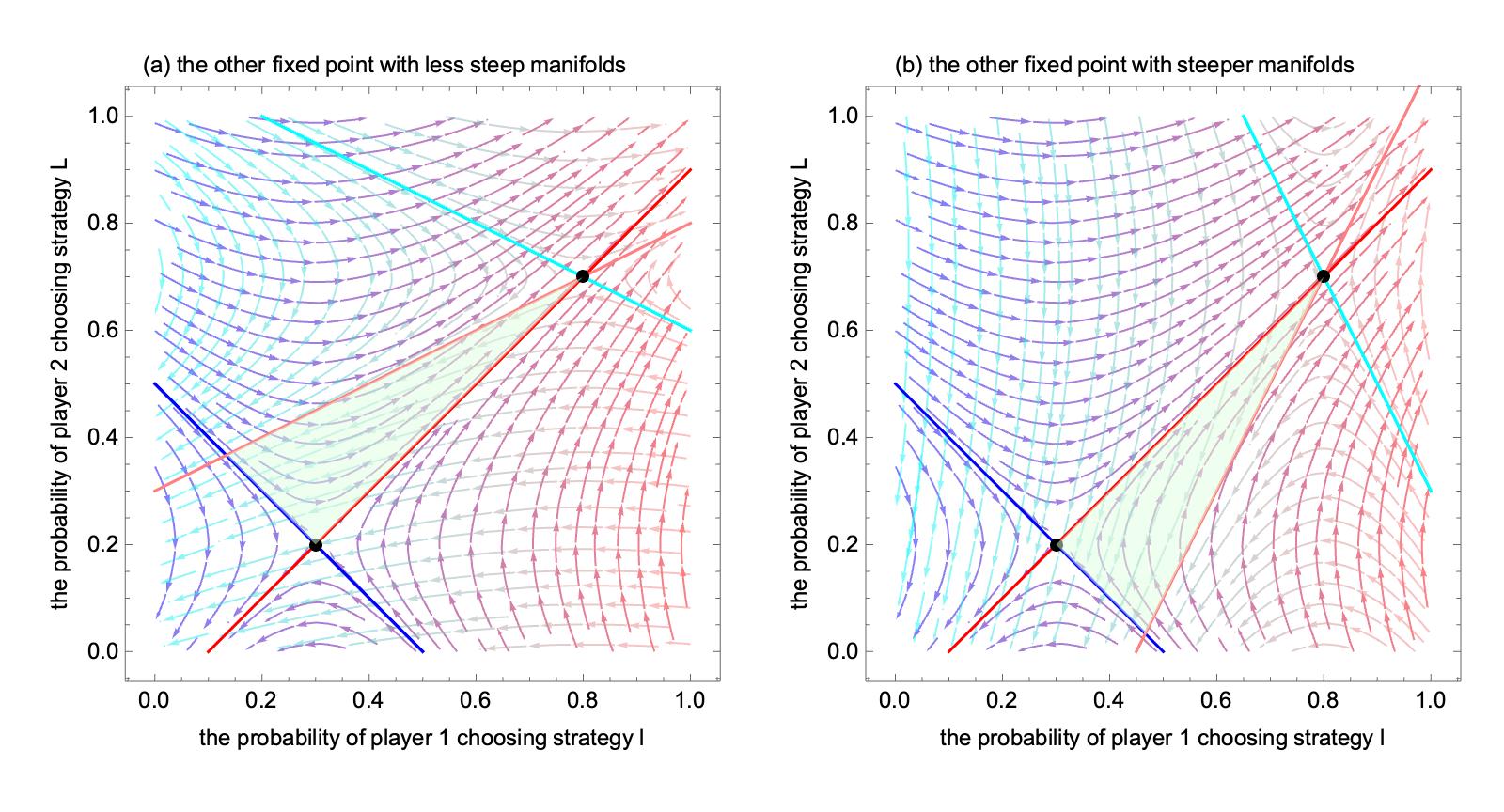}
    \caption{Linear Trapping: one fixed point on the unstable manifold of the other}
    \label{unstable-manifold}
\end{figure}

\begin{figure}[H]
    \centering
    \includegraphics[scale=0.5]{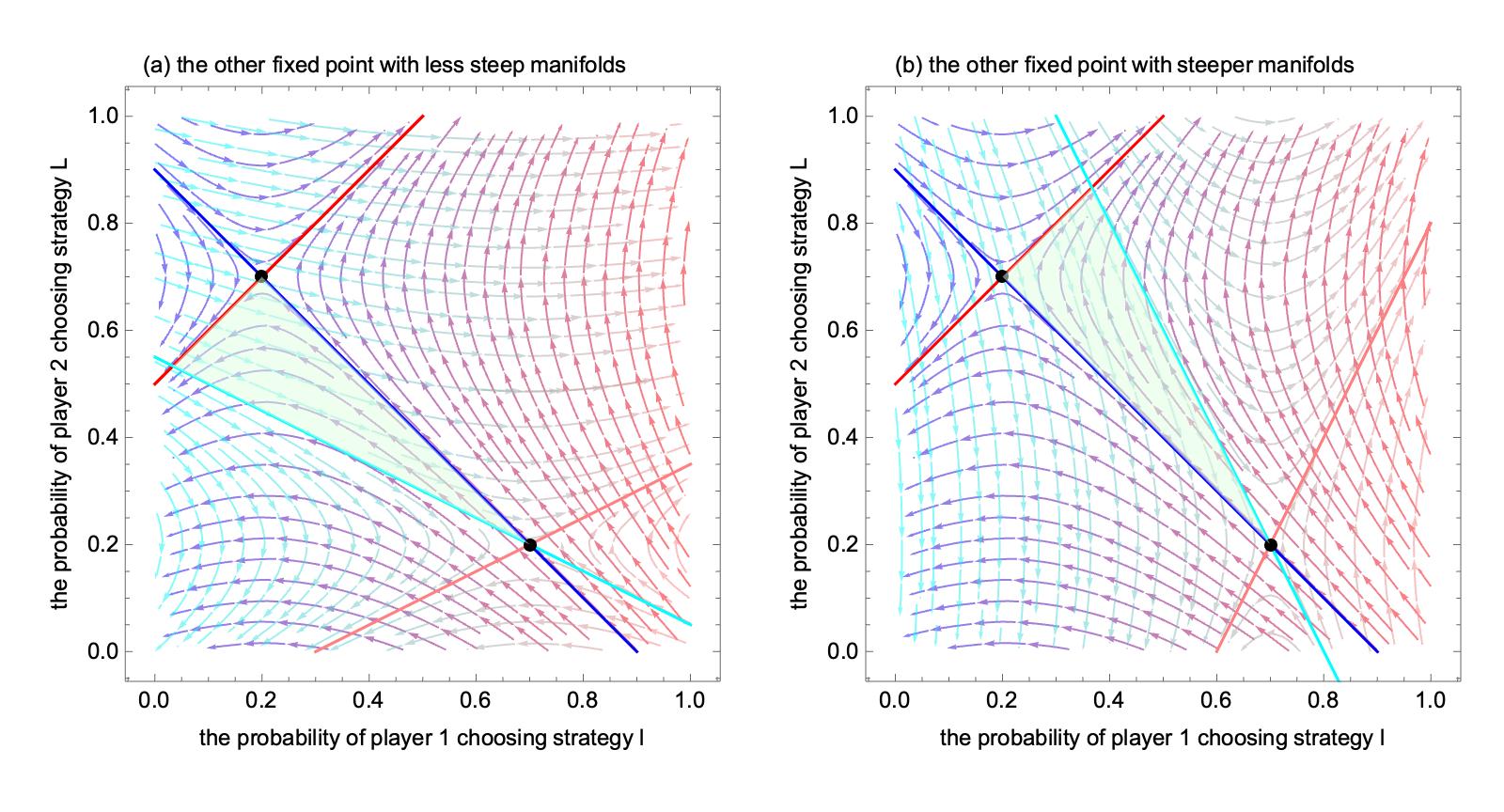}
    \caption{Linear Trapping: one fixed point on the stable manifold of the other}
    \label{stable-manifold}
\end{figure}

Still considering one fixed point on the stable or unstable manifold of the other fixed point, we now focus on the cases where their other manifolds are not parallel to each other. In Fig. \ref{unstable-manifold},  one fixed point is on the unstable manifold of the other fixed point, depending on whether the manifolds of the fixed point are less steep or steeper, we have two slightly different triangular trapping regions. And in Fig. \ref{stable-manifold}, one fixed point is on the stable manifold of the other fixed point. We end up with similar triangular trapping regions.

One might ask, what happens when

\begin{align}
   \min\{\sqrt{\frac{\beta}{\alpha}},\sqrt{\frac{\beta'}{\alpha'}}\}<\left|\frac{a^*-c^*}{b^*-d^*}\right|<\max\{\sqrt{\frac{\beta}{\alpha}},\sqrt{\frac{\beta'}{\alpha'}}\}?
   \label{extra}
\end{align}

Notice that \eqref{LR},\eqref{UD},\eqref{SUM} and \eqref{extra} include all possible configurations of the two distinct fixed points $E_1$ and $E_2$. In the case of \eqref{extra} we still have trapping regions but they are a little more complicated than the convex trapping regions we have shown so far. 

\begin{figure}[H]
    \centering
    \includegraphics[scale=0.5]{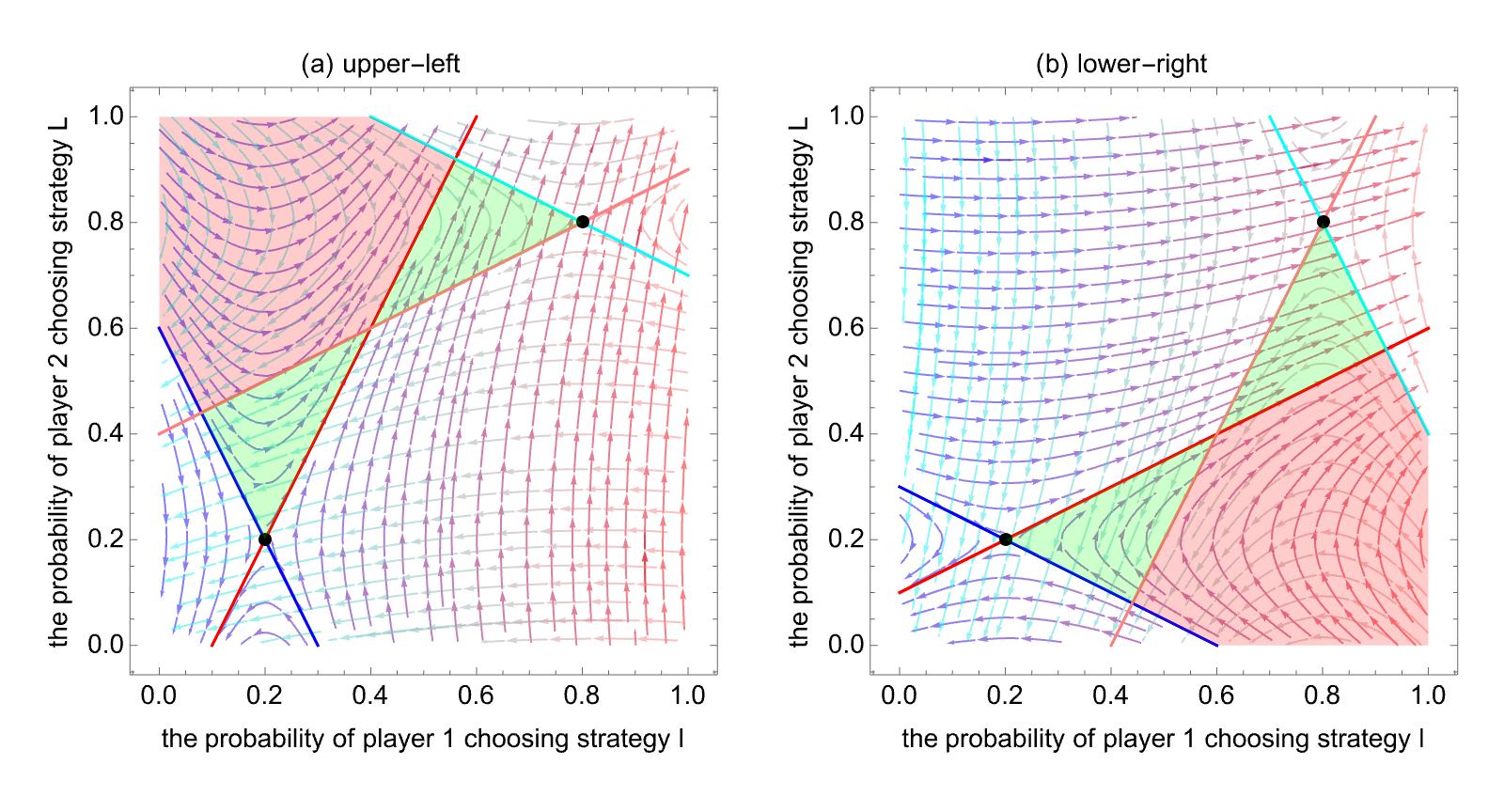}
    \caption{Linear Trapping: mixed configuration for fixed points positions}
    \label{mixed}
\end{figure}

Whenever we have the slope of the segment $E_1E_2$ smaller than the eigenvector slope of $E_1$, that means $E_2$ is either on the left or the right of the 4 parts of the phase space divided by the contracting and expanding directions of $E_1$. This is why \eqref{LR} corresponds to the case when $E_1$ and $E_2$ are in the "left-right" configuration for each other. Similarly \eqref{UD} corresponds to when the two fixed points are in the "up-down" configuration for each other. Intuitively \eqref{extra} means $E_1$ is in the "left-right" configuration w.r.t. $E_2$ while $E_2$ is in the "up-down" configuration w.r.t. $E_1$, or the other way around. We draw two different examples in Fig. \ref{mixed}.

Here we have a slightly more complicated mode of trapping. For example, in Fig. \ref{mixed}(a) each of the green triangular regions on their own is a trapping region. But the red shaded region combined with the two green triangles is also a trapping region. When we consider the combined shape, which is not convex any more, we could end up with a trapped loop that forms a butterfly shape or "8" shape. The trajectory starts in the red region, shoots into one of the green ears and then returns to the red region. If we switch the governing equations in time, then the trajectory is going to shoot into the other green ear and return to the red region. Although Fig. \ref{mixed}(a) and (b) have different locations for the trapping regions, the dynamics is the same. The phase portraits for cases where the red region is located at the upper-right corner and lower-left corner can also be drawn but we leave that up for the readers.

For the sake of completeness, we also put a simplified sketch of all possible linear trapping regions in Fig. \ref{list} below. We only draw the upper-left location trapping for the mixed configuration.

\begin{figure}[H]
    \centering
    \begin{tikzpicture}
    \draw[gray] (0,0) rectangle (4,4);
    \filldraw[black] (0.8,2.4) circle (1.5pt) (3.2,2.4) circle (1.5pt);
    \draw[red] (0,1.6) -- (2.4,4);
    \draw[blue] (0,3.2) -- (3.2,0);
    \draw[pink] (0,0.8) -- (4, 2.8);
    \draw[cyan] (0,4) -- (4,2);
    \draw[->] (1.6,1.8) .. controls (1,2.4) .. (1.6,3);
    \draw[->] (1.8,3) .. controls (3,2.4) .. (1.8,1.8);
    \draw[gray] (4.5,0) rectangle (8.5,4);
    \filldraw[black] (6.5,0.8) circle (1.5pt) (6.5,3.2) circle (1.5pt);
    \draw[red] (5.7,0) -- (8.5,2.8);
    \draw[blue] (4.5,2.8) -- (7.3,0);
    \draw[pink] (4.5,2.2) -- (8.1,4);
    \draw[cyan] (4.9,4) -- (8.5,2.2);
    \draw[->] (5.2,2.3) .. controls (6.5,1) .. (7.8,2.3);
    \draw[->] (7.8,2.4) .. controls (6.5,3) .. (5.2,2.4);
    \draw[gray] (9,0) rectangle (13,4);
    \filldraw[black] (10,3) circle (1.5pt) (12,1) circle (1.5pt);
    \draw[blue] (9,4) -- (13,0);
    \draw[red] (9,2) -- (11,4);
    \draw[red] (11,0) -- (13,2);
    \draw[->] (10.5,2.6) -- (11.5,1.6);
    \draw[->] (11.5,1.4) -- (10.5,2.4);
    \draw[gray] (0,-4.5) rectangle (4,-0.5);
    \filldraw[black] (1,-4) circle (1.5pt) (3,-2) circle (1.5pt);
    \draw[red] (0.5,-4.5) -- (4,-1);
    \draw[blue] (0,-3) -- (1.5,-4.5);
    \draw[blue] (1.5,-0.5) -- (4,-3);
    \draw[->] (1.5,-3.4) -- (2.5,-2.4);
    \draw[->] (2.5,-2.6) -- (1.5,-3.6);
    \draw[gray] (4.5,-4.5) rectangle (8.5,-0.5);
    \filldraw[black] (6,-3.5) circle (1.5pt) (7.5,-2) circle (1.5pt);
    \draw[blue] (4.5,-2) -- (7,-4.5);
    \draw[red] (5,-4.5) -- (8.5,-1);
    \draw[pink] (4.5, -3.5) -- (8.5,-1.5);
    \draw[cyan] (4.5,-0.5) -- (8.5,-2.5);
    \draw[->] (5.6,-3) .. controls (6,-3.4) .. (7,-2.4);
    \draw[->] (6.7,-2.5) -- (5.9,-2.9);
    \draw[gray] (9,-4.5) rectangle (13,-0.5);
    \filldraw[black] (10.5,-3.5) circle (1.5pt) (12,-2) circle (1.5pt);
    \draw[blue] (9,-2) -- (11.5,-4.5);
    \draw[red] (9.5,-4.5) -- (13,-1);
    \draw[pink] (10.75,-4.5) -- (12.75,-0.5);
    \draw[cyan] (11.25,-0.5) -- (13,-4);
    \draw[->] (10.9,-3.8) .. controls (10.6,-3.5) .. (11.6,-2.5);
    \draw[->] (11.45,-2.9) -- (11.05,-3.7);
    \draw[gray] (0,-9) rectangle (4,-5);
    \filldraw[black] (1,-6.5) circle (1.5pt) (2.5,-8) circle (1.5pt);
    \draw[blue] (0,-5.5) -- (3.5,-9);
    \draw[red] (0,-7.5) -- (2.5,-5);
    \draw[cyan] (0,-6.75) -- (4,-8.75);
    \draw[pink] (0.5,-9) -- (4,-7.25);
    \draw[->] (2,-7.6) .. controls (1,-6.6) .. (0.6,-7);
    \draw[->] (0.9,-7.1) -- (1.7,-7.5);
    \draw[gray] (4.5,-9) rectangle (8.5,-5);
    \filldraw[black] (5.5,-6.5) circle (1.5pt) (7,-8) circle (1.5pt);
    \draw[blue] (4.5,-5.5) -- (8,-9);
    \draw[red] (4.5,-7.5) -- (7,-5);
    \draw[cyan] (5.5,-5) -- (7.5,-9);
    \draw[pink] (6.5,-9) -- (8.5,-5);
    \draw[->] (6.5,-7.4) .. controls (5.6,-6.5) .. (6,-6.1);
    \draw[->] (6.1,-6.4) -- (6.4,-7);
    \draw[gray] (9,-9) rectangle (13,-5);
    \filldraw[black] (9.8,-8.2) circle (1.5pt) (12.2,-5.8) circle (1.5pt);
    \draw[blue] (9,-6.6) -- (10.2,-9);
    \draw[red] (9.4,-9) -- (11.4,-5);
    \draw[cyan] (10.6,-5) -- (13,-6.2);
    \draw[pink] (9,-7.4) -- (13,-5.4);
    \draw[->] (11,-5.3) .. controls (12.3,-5.8) .. (9.5,-7.05);
    \draw[->] (9.3,-7) .. controls (9.8,-8.3) .. (11.1,-5.4);
    \draw[->] (11.2,-5.5) .. controls (12,-5.8) .. (10.8,-6.4);
    \draw[->] (10.85,-6.3) -- (11.2,-5.6);
    \draw[->] (9.5,-7.3) .. controls (9.8,-8) .. (10.4,-6.8);
    \draw[->] (10.3,-6.85) -- (9.6,-7.2);
    \end{tikzpicture}
    \caption{All cases for linear trapping}
    \label{list}
\end{figure}
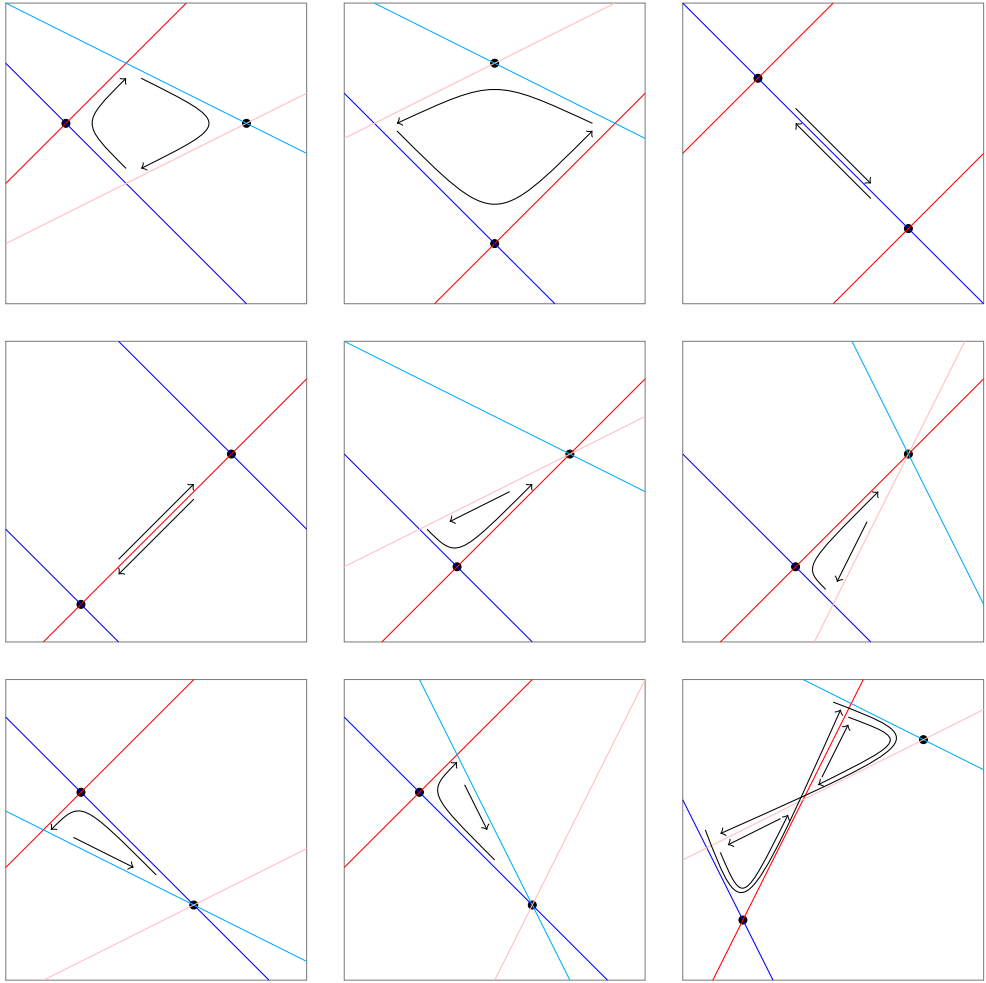

\section{Constant of motion for the system}

Our system is governed by these equations

$$
\begin{cases}
\dot x_1=x_1(1-x_1)\{[a_{11}(t)+a_{22}(t)-a_{12}(t)-a_{21}(t)]y_1-[a_{22}(t)-a_{12}(t)]\},\\
\dot y_1=y_1(1-y_1)\{[b_{11}(t)+b_{22}(t)-b_{12}(t)-b_{21}(t)]x_1-[b_{22}(t)-b_{21}(t)]\}.
\end{cases}
$$

For convenience, let's assume $a_{ij},b_{ij}$'s are constants here, for $i,j=1,2$. Let $p=a_{11}+a_{22}-a_{12}-a_{21}$, $q=a_{22}-a_{12}$, $u=b_{11}+b_{22}-b_{12}-b_{21}$, $v=b_{22}-b_{21}$. We will write $x$ for $x_1$, $y$ for $y_1$ in the following solving process to make the symbols easier to read and write.

Our system becomes

$$
\dot x=x(1-x)(py-q),
$$
$$\dot y=y(1-y)(ux-v).$$

Divide the second equation by the first, we get

\begin{align*}
    \frac{dy}{dx}=\frac{y(1-y)(ux-v)}{x(1-x)(py-q)}.
\end{align*}

This is a separable equation.

\begin{align*}
    &\frac{(py-q)dy}{y(1-y)}=(\frac{p-q}{1-y}-\frac{q}{y})dy=\frac{(ux-v)dx}{x(1-x)}=(\frac{u-v}{1-x}-\frac{v}{x})dx,\\
    &(q-p)\ln(1-y)-q\ln y=(v-u)\ln(1-x)-v\ln x+c,\\
    &\frac{(1-y)^{q-p}}{y^q}=c\frac{(1-x)^{v-u}}{x^v}, x^v(1-y)^{q-p}=cy^q(1-x)^{v-u},\\
    & V(x,y)=x^v(1-x)^{u-v}y^{-q}(1-y)^{q-p} = c
\end{align*}

$V(x,y)=c$ is a constant of motion \cite{hofbauer} for our system. It gives the analytic implicit solution to our system. Without an explicit solution, we still don't have the exact behavior of $x(t)$ and $y(t)$.

\section{Discussion}

Understanding how cooperative traits are maintained in populations is an essential problem of deep importance not only in evolutionary biology, but also in microbial ecology and systems biology. Ecological factors such as population density, disturbance frequency, population dispersal, resource supply, spatial structuring of populations, the presence of mutualisms, or the presence of a competing species in the environment, often play an important role in the evolution of cooperation. The effect that these and other ecological factors play on the evolution of cooperation is in general well understood. However, the reverse process, i.e., the effect that the evolution of cooperative traits may have on the ecological properties of populations is worthy of investigation using the framework of eco-evolutionary dynamics.

In this paper we considered eco-evolutionary games that incorporate the impact of switching environment on game dynamics. We consider a changing environment, which is represented as a switched system between two different payoff matrices. Motivated by the study of the tragedy of the commons in evolutionary biology, we demonstrate how alternative dynamics can arise. By switching between two environments, we can achieve an oscillating dynamic equilibrium and maintain coexistence which is impossible otherwise. Cooperators and defectors can coexist if the system state does not exceed the trapping regions.

Our model uses simplifying assumptions and is just a first step toward a more sophisticated eco-evolutionary model~\cite{tilman2021evolution,jnawali2022stochasticity}. Still, we strongly believe that our results provide valuable insights in designing control mechanisms. We may consider relaxing or extending these simplifying assumptions for future research. For example, instead of a switch between two different payoff matrices, we could consider a linear combination of the two. The possibilities under a general nonlinear payoff matrix are likely even broader. We could also consider incorporating multiplayer games into our current model toward a more realistic description of population dynamics.

Other possible extensions of our work could be to consider feedback between strategies and patch quality~\cite{fahimipour2022sharp,liu2022environmental}, where not all players experience the same environment. Negative environmental feedback gives immediate benefits in rich environments but results in net costs in poor environments.

In the analysis of the trapping regions we find multiple segments of stable and unstable manifolds in phase graphs, which is a new dynamical phenomenon that is different from interior fixed equilibrium situation obtained in previous works. We provide constructive specific control examples where the system states evolve contained inside the trapping regions. In some very special cases we even get an equilibrium manifold, depending on the initial condition, our system could be trapped on the equilibrium manifold. This combination of eco-evolutionary dynamics and control theory result could be used to help steer the population to desired states in social dilemma games. Although our present study is focused on mixed strategy bimatrix games, system properties are unchanged in a small neighbourhood, our results on control can be generalized to many other important situations.




\section*{Acknowledgements}
F.F. is supported by the Bill \& Melinda Gates Foundation (award no. OPP1217336), the NIH COBRE Program (grant no. 1P20GM130454), a Neukom CompX Faculty Grant, the Dartmouth Faculty Startup Fund, and the Walter \& Constance Burke Research Initiation Awar d.

\bibliographystyle{elsarticle-num} 

\begin{thebibliography}{100}

\bibitem{hauert2005game}
C. Hauret, G. Szab{\'o}, Game theory and physics, American Journal of Physics 73 (5) (2005) 405-414.

\bibitem{nowak2006evolutionary}
M. A. Nowak, Evolutionary dynamics: exploring the equations of life, Harvard University Press, 2006.

\bibitem{roca2009evolutionary}
C. P. Roca, J. A. Cuesta, A. S{\'a}nchez, Evolutionary game theory: Temporal and spatial effects beyond replicator dynamics, Physics of life reviews 6 (4) (2009) 208-249.

\bibitem{cressman2014replicator}
R. Cressman, Y. Tao, The replicator equation and other game dynamics, Proceedings of the National Academy of Sciences 111 (supplement 3) (2014) 10810-10817.

\bibitem{perc2010coevolutionary}
M. Perc, A. Szolnoki, Coevolutionary games-a mini review, BioSystems 99 (2) (2010) 109-125.

\bibitem{sanchez2013feedback}
A. S{\'a}nchez, J. Gore, Feedback between population and evolutionary dynamics determines the fate of social microbial populations, PLoS biology 11 (4) (2013) e1001547.

\bibitem{wang2020steering}
X. Wang, Z. Zhang, F. Fu, Steering eco-evolutionary game dynamics with manifold control, Proceedings of the Royal Society A 476 (2233) (2020) 20190643.

\bibitem{estrela2019environmentally}
S. Estrela, E. Libby, J. Van Cleve, F. D{\'e}barre, M. Deforet, W. R. Harcombe, J. Pe{\~n}a, S. P. Brown, M. E. Hochberg, Environmentally mediated social dilemmas, Trends in ecology \& evolution 34 (1) (2019) 6-18.

\bibitem{rand2021geometry}
D. A. Rand, A. Raju, M. Saez, F. Corson, E. D. Siggia, Geometry of gene regulatory dynamics, arXiv preprint arXiv:2105.13722 (2021).

\bibitem{weitz2016oscillating}
J. S. Weitz, C. Eksin, K. Paarporn, S. P. Brown, W. C. Ratcliff, An oscillating tragedy of the commons in replicator dynamics with game-environment feedback, Proceedings of the National Academy of Sciences 113 (47) (2016) E7518-E7525.

\bibitem{hilbe2018evolution}
C. Hilbe, {\v{S}}. {\v{S}}imsa, K. Chatterjee, M. A. Nowak, Evolution of cooperation in stochastic games, Nature 559 (7713) (2018) 246-249.

\bibitem{chen2018punishment}
X. Chen, A. Szolnoki, Punishment and inspection for governing the commons in a feedback-evolving game, PLoS computational biology 14 (7) (2018) e1006347.

\bibitem{shao2019evolutionary}
Y. Shao, X. Wang, F. Fu, Evolutionary dynamics of group cooperation with asymmetrical environmental feedback, EPL (Europhysics Letters) 126 (4) (2019) 40005.

\bibitem{hauert2019asymmetric}
C. Hauert, C. Saade, A. McAvoy, Asymmetric evolutionary games with environmental feedback, Journal of theoretical biology 462 (2019) 347-360.

\bibitem{tilman2020evolutionary}
A. R. Tilman, J. B. Plotkin, E. Ak{\c{c}}ay, Evolutionary games with environmental feedbacks, Nature communications 11 (1) (2020) 1-11.

\bibitem{hardin1968tragedy}
G. Hardin, The tragedy of the commons: the population problem has no technical solution; it requires a fundamental extension in morality, Science 162 (3859) (1968) 1243-1248.

\bibitem{wu2018coevolutionary}
T. Wu, F. Fu, L. Wang, Coevolutionary dynamics of aspiration and strategy in spatial repeated public goods games, New Journal of Physics 20 (6) (2018) 063007.

\bibitem{perez2022cooperation}
H. Perez, C. Gracia-Lazaro, F. Dercole, Y. Moreno, Cooperation in costly-access environments, New Journal of Physics (2022).

\bibitem{wang2022replicator}
X. Wang, M. Perc, Replicator dynamics of public goods games with global exclusion, Chaos: An Interdisciplinary Journal of Nonlinear Science 32 (7) (2022) 073132.

\bibitem{wang2022decentralized}
S. Wang, X. Chen, Z. Xiao, A. Szolnoki, Decentralized incentives for general well-being in networked public goods game, Applied Mathematics and Computation 431 (2022) 127308.

\bibitem{tarnita2017ecology}
C. E. Tarnita, The ecology and evolution of social behavior in microbes, Journal of Experimental Biology 220 (1) (2017) 18-24.

\bibitem{cardillo2020critical}
A. Cardillo, N. Masuda, Critical mass effect in evolutionary games triggered by zealots, Physical Review Research 2 (2) (2020) 023305.

\bibitem{chen2019imperfect}
X. Chen, F. Fu, Imperfect vaccine and hesteresis, Proceedings of the royal society B 286 (1894) (2019) 20182406.

\bibitem{ostrom1993coping}
E. Ostrom, R. Gardner, Coping with asymmetries in the commons: self-governing irrigation systems can work, Journal of economic perspectives 7 (4) (1993) 93-112.

\bibitem{chen2018social}
X. Chen, F. Fu, Social learning of prescribing behavior can promote population optimum of antibiotic use, Frontiers in Physics 6 (2018) 139.

\bibitem{cao2021eco}
L. Cao, B. Wu, Eco-evolutionary dynamics with payoff-dependent environmental feedback, Chaos, Solitons \& Fractals 150 (2021) 111088.

\bibitem{glaubitz2020oscillatory}
A. Glaubitz, F. Fu, Oscillatory dynamics in the dilemma of social distancing, Proceedings of the Royal Society A 476 (2243) (2020) 20200686.

\bibitem{hofbauer}
J. Hofbauer, K. Sigmund, Evolutionary Games and Population Dynamics, Cambridge University Press, 1998.

\bibitem{tilman2021evolution}
A. R. Tilman, V. V. Vasconcelos, E. Akcay, J. B. Plotkin, The evolution of forecasting for decision making in dynamic environments, arXiv preprint arXiv:2108.00047 (2021).

\bibitem{jnawali2022stochasticity}
K. Jnawali, M. Anand, C. T. Bauch, Stochasticity-induced persistence in coupled social-ecological systems, Journal of Theoretical Biology 542 (2022) 111088.

\bibitem{fahimipour2022sharp}
A. K. Fahimipour, F. Zeng, M. Homer, A. Traulsen, S. A. Levin, T. Gross, Sharp thresholds limit the benefit of defector avoidance in cooperation on networks, Proceedings of the National Academy of Sciences 119 (33) (2022) e2120120119.

\bibitem{liu2022environmental}
F. Liu, B. Wu, Environmental quality and population welfare in markovian eco-evolutionary dynamics, Applied Mathematics and Computation 431 (2022) 127309.

\end{thebibliography}


\end{document}